\documentclass[12pt]{amsart}
\usepackage{amsmath}
\usepackage{amscd}
\usepackage[psamsfonts]{amssymb}
\usepackage[mathscr,psamsfonts]{eucal}
\DeclareMathAlphabet{\eurm}{U}{eur}{m}{n}
\newcommand{\F}{{\mathcal F}(E_1,E_2)}
\newcommand{\mc}[1]{{\mathcal{#1}}}
\DeclareMathOperator{\id}{{id}}
\DeclareMathOperator{\Aut}{{Aut}}
\DeclareMathOperator{\byd}{\,{\raisebox{.1ex}{$\eurm :$}{\eurm =}}\,}
\newcommand{\sig}{\sigma}
\newcommand{\alp}{\alpha}
\newcommand{\gam}{\gamma}
\newcommand{\bEq}{\begin{eqnarray}}
\newcommand{\eEq}{\end{eqnarray}}
\newcommand{\beq}{\begin{eqnarray*}}
\newcommand{\eeq}{\end{eqnarray*}}

\newcommand{\com}{\circ}
\newcommand{\car}{\times}
\newcommand{\ten}{\otimes}
\newcommand{\Rn}{{\mathbb R}}
\newcommand{\END}{{\,\text{\qedsymbol}}}
\newcommand{\mto}{\mapsto}

\begin{document}

\newcounter{theorem}

\newtheorem{definition}[theorem]{Definition}
\newtheorem{lemma}[theorem]{Lemma}
\newtheorem{proposition}[theorem]{Proposition}
\newtheorem{theorem}[theorem]{Theorem}
\newtheorem{corollary}[theorem]{Corollary}
\newtheorem{remark}[theorem]{Remark}
\newtheorem{example}[theorem]{Example}
\newtheorem{Note}[theorem]{Note}
\newcounter{assump}
\newtheorem{Assumption}{\indent Assumption}[assump]
\renewcommand{\thetheorem}{\thesection.\arabic{theorem}}

\newcommand{\bCr}{\begin{corollary}}
\newcommand{\eCr}{\end{corollary}}
\newcommand{\bDf}{\begin{definition}\em}
\newcommand{\eDf}{\end{definition}}
\newcommand{\bLm}{\begin{lemma}}
\newcommand{\eLm}{\end{lemma}}
\newcommand{\bPr}{\begin{proposition}}
\newcommand{\ePr}{\end{proposition}}
\newcommand{\bRm}{\begin{remark}\em}
\newcommand{\eRm}{\end{remark}}
\newcommand{\bEx}{\begin{example}\em}
\newcommand{\eEx}{\end{example}}
\newcommand{\bTh}{\begin{theorem}}
\newcommand{\eTh}{\end{theorem}}
\newcommand{\bNt}{\begin{Note}\em}
\newcommand{\eNt}{\end{Note}}
\newcommand{\bPf}{\begin{proof}[\noindent\indent{\sc Proof}]}
\newcommand{\ePf}{\renewcommand{\qedsymbol}{}\end{proof}}
\title{Functorial prolongations of some functional bundles}

\author[A. Cabras, J. Jany\v ska, I. Kol\'a\v{r}]
        {Antonella Cabras, Josef Jany\v ska, Ivan Kol\'a\v{r}
}

\dedicatory{To Andrzej Zajtz, on the occasion of his 70th birthday.}



\address{
\newline
Antonella Cabrass
\newline
Department of Applied Mathematics, Florence University
\newline
Via S. Marta 3, 50139 Florence, Italy
\newline
E-mail: {\tt cabras@dma.unifi.it}
\newline
{\ }
\newline
Josef Jany\v{s}ka
\newline
Department of Mathematics, Masaryk University
\newline
Jan\'a\v ckovo n\'am. 2a, 662 95 Brno, Czech Republic
\newline
E-mail: {\tt janyska@math.muni.cz}
\newline
{\ }
\newline
Ivan Kol\'a\v{r}
\newline
Department of Algebra and Geometry, Masaryk University
\newline
Jan\'a\v ckovo n\'am. 2a, 662 95 Brno, Czech Republic
\newline
E-mail: {\tt kolar@math.muni.cz}
}

\thanks{The second and third authors were supported
by the Ministry of Education of the Czech Republic under
the project
MSM 143100009.}

\begin{abstract}
We discuss two kinds of functorial prolongations of the functional
bundle of all smooth maps between the fibers over the same base point
of two fibered manifolds over the same base. We study
the prolongation of vector fields in both cases and we prove that
the bracket is preserved. Our proof is based on several new
results concerning the finite dimensional Weil bundles.
\end{abstract}
\maketitle

\section*{Introduction}
\setcounter{equation}{0}
\setcounter{theorem}{0}
Let $E_1$ and $E_2$ be two classical fiber bundles over the same
base $M$. The differential geometric investigation of the functional
bundle $\F\to M$ of all smooth maps from a fiber of $E_1$ into
the fiber of $E_2$ over the same base point was iniciated by the paper
by A. Jadczyk and M. Modugno on the Schr\"odinger connection,
\cite{JadMod92}, \cite{JanMod02}. The simpliest cases
of the tangent bundle $T\F\to TM$ and of the $r$-th jet prolongation
$J^r\F\to M$ are discussed in \cite{CabKol95}. In the present paper
we first clarify that the essential assumption for these constructions
is that $T$ is a product preserving bundle functor on the classical
category
$\mc{M}f$ of all smooth manifolds and all smooth maps and $J^r$
is a fiber product preserving bundle functor on the category
${\mc{F}}{\mc{M}}_m$ of all fibered manifolds
with $m$-dimensional bases
and of all fibered manifold morphisms covering local diffeomorphisms.
Every product preserving bundle functor $F$ on $\mc{M}f$ is a
Weil functor $F= T^A$, where $A$ is a Weil algebra,
\cite{KolMicSlo93}. The general construction of $T^A\F\to T^A M$
was presented by the third author in \cite{Kol96}, \cite{Kol96a},
see also Section 2 of the present paper. We underline that this
construction is based on the covariant approach to Weil bundles
and their natural transformations, \cite{Kol86}, \cite{KolMicSlo93}.
On the other hand, in \cite{KolMik99} it was deduced that every
fiber product preserving bundle functor $G$ on
${\mc{F}}{\mc{M}}_m$
is of the form $G=(A,H,t)$, where $A$ is a Weil algebra, $H$
is a group homomorphism $H:G^r_m\to \Aut A$ of the $r$-th jet
group $G^r_m$ in dimension $m$ into the group of all algebra
automorphisms of $A$ and $t:{\mathbb D\,}{}^r_m\to A$ is an equivariant
algebra homomorphism, where
${\mathbb D\,}{}^r_m= J^r_0({\mathbb R}^m,{\mathbb R})$
is the Weil algebra corresponding to the functor of
$(m,r)$-velocities. In Section 6 of the present paper
we construct $G\F\to M$ in a way that generalizes the case
of $J^r\F\to M$.

Our main geometric problem is the prolongation of vector fields
on $\F$ with respect to $F$ and $G$. Since we cannot use the flow
in the functional case, we start from the fact that the classical
flow prolongation with respect to $T^A$ of a vector field $M\to TM$
coincides with the composition of its $T^A$-prolongation
$T^AM\to T^ATM$ with the exchange map $\kappa^A_M:T^ATM\to TT^AM$.
We apply this idea to a vector field $X$ on $\F$ and we
say the composition ${\mc{T}}^AX=\kappa^A_{\F} \com T^AX$
to be the field prolongation of $X$. The bracket of vector fields
on $\F$ is defined in terms of the strong difference,
\cite{CabKol95}, \cite{KolMicSlo93}. Proposition \ref{Pr1}
in Section 3 reads
that $\mc{T}^A$ preserves the bracket of vector fields even in the
functional case. To deduce it, we develop, in Sections 4 and 5,
a purely algebraic proof of the fact that ${\mc{T}}^A$
preserves bracket
in the manifold case. For this purpose we need certain new lemmas
concerning the classical Weil bundles, which are collected
in Sections 4 and 5. In particular, we present a complete
description of the strong difference in terms of Weil algebras.
In Section 7 we study the prolongation of vector fields to $G\F$
and we prove that the bracket is preserved even in this case.
Finally we remark that  an  interesting kind of exchange morphism,
which was introduced recently for the manifold case
in \cite{Kol03}, can
be extended to the functional bundles as well.

In Section 1 we present a simplified version of the theory of
smooth spaces in the sense of A. Fr\"olicher, \cite{Fro82}, which
we call $F$-smooth spaces, and of $F$-smooth bundles.
Special attention is paid to the functorial character of the
construction of $\F$ and to the concept of finite order morphism.

If we deal with finite dimensional manifolds and maps between them,
we always assume they are of class $C^\infty$, i.e. smooth
in the classical sense. Unless otherwise specified,
we use the terminology and notation from the monograph
\cite{KolMicSlo93}.

\section{$F$-smooth bundles}
\setcounter{equation}{0}
\setcounter{theorem}{0}
We shall use the following simplified version, \cite{CabKol98}, of the
theory of smooth spaces by A. Fr\"olicher, \cite{Fro82}.

\bDf
An $F$-smooth space is a set $S$ along with a set $C_{S}$
of maps $c:\Rn\to S$, which are called $F$-smooth curves,
satisfying the following two conditions:

(i) each constant curve $\Rn \to S$ belongs to $C_{S}$,

(ii) if $c\in C_{S}$ and $\gam\in C^\infty(\Rn,\Rn)$,
then $c\com \gam\in C_{S}$.

\noindent
If $(S',C_{S'})$ is another $F$-smooth space, a map
$f:S\to S'$ is said to be $F$-smooth, if $f\com c$ is an
$F$-smooth curve on $S'$ for every $F$-smooth curve $c$ on $S$.
\END
\eDf

So we obtain the category ${\mc{S}}$ of $F$-smooth spaces. Every
subset $\bar{S}\subset S$ is also an $F$-smooth space,
if we define $C_{\bar{S}}\subset C_{S}$ to be the subset of
the curves with values in $\bar{S}$. In particular every smooth
manifold $M$ turns out to be an $F$-smooth space by assuming
as $F$-smooth curves just the smooth curves. Moreover, a map between
smooth manifolds is $F$-smooth, if and only if it is smooth.

We find it useful to define the concept of $F$-smooth bundle in
a more general form than in \cite{CabKol98}.

\bDf
An $F$-smooth bundle is a triple of an $F$-smooth space $S$,
a smooth manifold $M$ and a surjective $F$-smooth map
$p:S\to M$. If $p':S'\to M'$ is another
$F$-smooth bundle, then a morphism of $S$ into $S'$ is a pair of
an $F$-smooth map $f:S\to S'$ and a smooth map $\underline f:
M\to M'$ satisfying $\underline f\com p= p'\com f$.
\END
\eDf

Thus we obtain the category ${\mc{S}\mc{B}}$ of
$F$-smooth bundles.
Every subset $\bar{S}\subset S$ satisfying $p(\bar{S})=M$
is also an  $F$-smooth bundle.

An important class of  $F$-smooth bundles are the bundles
of smooth maps between the fibers over the same base point
of two classical
fibered manifolds $p_1: E_1\to M$ and  $p_2: E_2\to M$.
We write
\beq
\F=\underset{x\in M}{\bigcup}
        C^\infty (E_{1x},E_{2x})
\eeq
and denote by $p:\F\to M$ the canonical
projection. A curve $c:\Rn\to \F$ is called
$F$-smooth, if $\underline{c}\byd p\com c:\Rn\to M$ is a smooth
map and the induced map
\beq
\tilde{c}:\underline{c}^* E_1\to E_2,\quad \tilde{c}(t,y)=
        c(t)(y),\quad p_1(y)=\underline{c}(t)
\eeq
is also smooth, \cite{CabKol95}.

Write ${\mc{F}\mc{M}}^I\subset{\mc{F}\mc{M}}$
for the subcategory of locally trivial fibered manifolds
whose morphisms are diffeomorphisms on the fibers.
Let ${\mc{F}\mc{M}}^I\times_{\mc{B}} {\mc{F}\mc{M}}$ denote
the category whose objects are pairs $(E_1,E_2)$ with $E_1\to M$ in
${\mc{F}\mc{M}}^I$ and  $E_2\to M$ in
${\mc{F}\mc{M}}$ and morphisms are pairs $(f_1,f_2)$ with
$f_1:E_1\to E_3$ in
${\mc{F}\mc{M}}^I$ and  $f_2:E_2\to E_4$ in
${\mc{F}\mc{M}}$ over the same base map
$\underline{f}:M\to N$, where $N$ is the common base of $E_3$
and $E_4$. If we define ${\mc{F}}(f_1,f_2): \F \to
{\mc{F}}(E_3,E_4)$ by
\bEq
{\mc{F}}(f_1,f_2)(h)= f_2(x)\com h\com f_1^{-1} (\underline{f}(x)),
        \quad h\in C^\infty(E_{1x}, E_{2x}),
\eEq
then $\mc{F}$ is a functor on
${\mc{F}\mc{M}}^I\times_{\mc{B}} {\mc{F}\mc{M}}$
with values in the category ${\mc{S}}{\mc{B}}$.

\bDf
Every $F$-smooth subbundle $S\subset \F$ will be called
a functional $F$-smooth bundle.
\END
\eDf

If $S'\subset {\mc{F}}(E_3,E_4)$ is another functional
$F$-smooth bundle and $(f_1,f_2)$ has the property
${\mc{F}}(f_1,f_2)(S)\subset S'$, then the restricted
and corestricted map will be interpreted as
an ${\mc{S}}{\mc{B}}$-morphism $S\to S'$.

Consider a smooth map $q:E_3\to E_1$.

\bDf
An ${\mc{S}}{\mc{B}}$-morphism $D:\F\to {\mc{F}}(E_3,E_4)$
is said to be of the order $r$, if for every
$\varphi, \psi:E_{1x}\to E_{2x}$ and $v\in E_3$, $p_1(q(v))=x$,
\bEq\label{e2}
j^r_{q(v)}\varphi=j^r_{q(v)}\psi\quad\text{implies}\quad
        D(\varphi)(v)= D(\psi)(v)\,.\quad\END
\eEq
\eDf

Consider the fibered manifold
\bEq\label{e3}
{\mc{F}}J^r(E_1,E_2)=  \underset{x\in M}{\bigcup}
        J^r(E_{1x}, E_{2x})\to E_1\,.
\eEq
By (\ref{e2}), $D$ induces the so called associated map
\beq
{\mc{D}}:{\mc{F}}J^r(E_1,E_2)\times_{E_1} E_3 \to E_4\,.
\eeq

\noindent
In the same way as in \cite{CabKol95} one proves that ${\mc{D}}$
is a smooth map.

We express the coordinate form of $\mc{D}$ in the case
$q:E_3\to E_1$ is an $\mc{F}\mc{M}$-morphism that is
a surjective submersion on each fiber of $E_3$.
Let $x^i$ or $u^a$ be some local coordinates on $M$ or $N$
and $y^p$ or $z^s$ or $(y^p,v^b)$ or $w^c$ be some additional
fiber coordinates on $E_1$ or $E_2$ or $E_3$ or $E_4$,
respectively. Then $z^s_\alp$ are the induced coordinates on
$\mc{F}J^r(E_1,E_2)$, where $0\le |\alp| \le r$ is
a multiindex, the range of which is the fiber dimension of $E_1$,
and the coordinate expression of $\mc{D}$ is
\bEq\label{e4}
u^a= f^a(x^i)\,,\quad w^c = f^c(x^i,y^p,z^s_\alp, v^b)\,,
\eEq
where $f^a$ and $f^c$ are smooth functions.

The concept of $r$-th order morphism can be modified to
a functional $F$-smooth bundle $S\subset \F$ analogously
to \cite{KolMicSlo93}, Section 18.

\section{The tangent-like case}
\setcounter{equation}{0}
\setcounter{theorem}{0}

Let $A$ be a Weil algebra of the width $k$. Under the covariant
approach, \cite{Kol86}, \cite{KolMicSlo93}, the elements
of a Weil bundle $T^AM$ are the $A$-velocities $j^Ag$ of smooth
maps $g:\Rn^k\to M$. For a smooth map $f:M\to N$, we define
$T^Af: T^AM\to T^AN$ by
\bEq\label{e5}
T^Af(j^Ag) = j^A(f\com g)\,.
\eEq
If $B$ is another Weil algebra of the width $l$, then every algebra
homomorphism $\mu:A\to B$ can be generated by a $B$-velocity $j^Bh$
of a map $h:\Rn^l\to \Rn^k$. The natural transformation
$\mu_M:T^AM\to T^BM$ induced by $\mu$ has the form of a reparametrization
\bEq\label{e6}
\mu_M(j^Ag) = j^B(g\com h)\,.
\eEq

Consider $\F$. We have $T^Ap_i:T^A E_i\to T^AM$ and we write
$T^A_X E_i\byd (T^Ap_i)^{-1}(X)$, $X\in T^AM$, $i=1,2$.
Let $g_1,g_2:\Rn^k\to \F$ be two $F$-smooth maps satisfying
$j^A(p\com g_1)= j^A(p\com g_2)\in T^AM$. Then we construct
the associated maps $T^A_0g_i: T^A_X E_1\to T^A_X E_2$,
\beq
T^A_0 g_i(j^Af(u))= j^A g_i(u)(f(u))\,,\quad u\in \Rn^k\,,
\eeq
where $f:\Rn^k\to E_1$ satisfies $p\com g_i=p_1\com f$, $i=1,2$.
If $T^A_0 g_1= T^A_0 g_2$, we say that $g_1$ and $g_2$ determine
the same $A$-velocity $j^Ag_1 = j^A g_2$. The set $T^A\F$ of all
such $A$-velicities is a subspace in $\mc{F}(T^A E_1,T^A E_2)\to T^AM$,
so a functional $F$-smooth bundle. In the product case
$E_i= M\times Q_i$, $i=1,2$, the third author deduced in \cite{Kol96}
\bEq\label{e7}
T^A(M\times Q_1, M\times Q_2) = T^A M\times C^\infty(Q_1, T^A Q_2)\,.
\eEq

In \cite{Kol96} it was also clarified that the idea
of reparametrization (\ref{e6}) can be applied
to $j^Ag\in T^A\F$ as well. So every algebra homomorphism
$\mu=j^Bh:A\to B$ induces an $F$-smooth map
\bEq\label{e8}
\mu_{\F} :T^A\F\to T^B\F\,,\quad j^Ag\mto j^B(g\com h)\,.
\eEq

Consider a functional $F$-smooth bundle $S\subset \F$. Then
$T^AS\subset T^A\F$ means the subset of all $j^Ag$, $g:\Rn^k\to S$.

\bDf
An $\mc{S}\mc{B}$-morphism $D:S\to \mc{F}(E_3,E_4)$ is called
$A$-differentiable, if the rule
\beq
T^AD(j^Ag)=j^A(D\com g)
\eeq
defines an $F$-smooth map $T^AS \to T^A\mc{F}(E_3,E_4)$.
We say $D$ is strongly differentiable, if it is $A$-differentiable
for every Weil algebra $A$.
\END\eDf

If $D$ is strongly differentiable, then $T^A D$ is also strongly
differentiable. Indeed, analogously to the finite dimensional
case one verifies easily $T^B(T^AD)=T^{B\ten A}D$. In particular,
every finite order morphism is strongly differentiable, for
its associated map is smooth. Further, each morphism
$\mc{F}(f_1,f_2)$ is strongly differentiable and we have
\beq
T^A\mc{F}(f_1,f_2)(j^Ag(u))=j^A\big(f_2(p(g(u)))\com g(u)
        \com f_1^{-1}(\underline{f}(p(g(u))))\big)\,.
\eeq

Thus, $T^A\mc{F}$ is a functor on the category
$\mc{F}\mc{M}^I\car_{\mc{B}}\mc{F}\mc{M}$ with values in
$\mc{S}\mc{B}$.

Analogously to the finite dimensional case, \cite{CabKol01},
we define an $A$-field on $\F$ as a strongly differentiable
section $\F\to T^A\F$. In the case $A=\mathbb D$ of the algebra
of dual numbers, we obtain a vector field $X:\F\to T\F$.

\section{Prolongation of vector fields}
\setcounter{equation}{0}
\setcounter{theorem}{0}

In the manifold case, the exchange algebra homomorphism
$\kappa^A:A\ten \mathbb D\to \mathbb D\ten A$ defines
a natural transformation $\kappa^A_M: T^ATM\to TT^AM$.
For a classical vector field $X:M\to TM$, its flow prolongation
$\mc{T}^AX:T^A M\to TT^AM$ coincides with $\kappa^A_M\com T^AX$,
\cite{KolMicSlo93}.
For a vector field $X:\F\to T\F$, we also can construct
$T^AX:T^A\F\to T^AT\F$ and apply $\kappa^A_{\F}:T^AT\F\to TT^A\F$.
In this way we obtain a vector field on $T^A\F$.

\bDf
The vector field $\mc{T}^AX\byd \kappa^A_{\F}\com T^AX$ will
be called the field prolongation of $X$.
\END\eDf

We recall that the bracket of two vector fields $X,Y$
on $\F$ was defined by using the strong difference, \cite{CabKol95},
\bEq\label{e9}
[X,Y]= ({T}Y\com X)\div  ({T}X\com Y)\,.
\eEq
(For classical vector fields $X,Y:M\to TM$, (\ref{e9}) coincides
with the classical bracket, \cite{CabKol95}.)
We are going to deduce

\bPr\label{Pr1}
For every vector fields $X,Y$ on $\F$,
\bEq\label{e10}
\mc{T}^A([X,Y])=[\mc{T}^AX,\mc{T}^A Y]\,.
\eEq
\ePr

The proof will be based on the algebraic results of
the next two sections.

\section{The algebraic form of the strong difference}
\setcounter{equation}{0}
\setcounter{theorem}{0}

Write $p^T_M: TM\to M$ for the bundle projection.
We recall that two elements $X,Y\in TT_xM$ satisfying
\bEq\label{e11}
p^T_{TM} X= Tp^T_M Y\,,\quad p^T_{TM} Y= Tp^T_M X
\eEq
determine the strong difference
\bEq\label{e12}
X\div  Y \in T_xM\,,
\eEq
\cite{KolMicSlo93}. Denote by $S M$ the domain of definition of the
strong difference, i.e. $S M\subset TTM\car_M TTM$ is the subset of
all pairs $(X,Y)$ satisfying (\ref{e11}), and by $\sig_M:S M\to TM$
the map (\ref{e12}). For every smooth map $f:M\to N$, one verifies
easily that $(TTf,TTf)$ transforms $S M$ into $S N$.
So we obtain a map
\beq
S f:S M\to S N
\eeq
and $S$ is a bundle functor on $\mc{M}f$. Moreover,
the strong difference map is a natural transformation
\bEq\label{e13}
\sig_M:S M\to TM\,.
\eEq

The fact $S \Rn^m= \overset{5}{\car}\Rn^m$ implies that
$S$ preserves products. Write $\mathbb{S}$ for the corresponding
Weil algebra. In general, the sum of two Weil algebras $A=\Rn\car N_A$
and $B=\Rn\car N_B$ is defined by
\beq
A+B= \Rn\car N_A\car N_B
\eeq
with the induced multiplication that satisfies $ab=0$
for all $a\in N_A$, $b\in N_B$. Clearly, we have
\beq
T^AM\car_MT^BM = T^{A+B} M\,.
\eeq

Write $\mathbb D=\{a_0+a_1e\}$, $e^2=0$. Then $TT$ corresponds to
$\mathbb{D}\ten\mathbb{D}$, which is linearly generated by
$1,e_1, e_2,e_1e_2$. Let $\{1, E_1,E_2,E_1 E_2\}$ be the linear
generators of another copy of $\mathbb{D}\ten\mathbb{D}$. So
$\mathbb{S}$ is a subalgebra of $\mathbb{D}\ten\mathbb{D}+
\mathbb{D}\ten\mathbb{D}$ and (\ref{e11}) implies directly that
the elements of $\mathbb{S}$ are of the form
\beq
X=a_0+a_1(e_1+E_2)+a_2(e_2+E_1)+a_3e_1e_2+a_4E_1 E_2\,,
\eeq
$a_0,\dots,a_4\in \Rn$. By the definition of the strong difference,
\cite{KolMicSlo93}, the algebra homomorphism
$\sig:\mathbb{S}\to \mathbb{D}$ corresponding
to (\ref{e12}) is
\bEq\label{e14}
\sig(X)= a_0 +(a_3-a_4)e\,.
\eEq

Write $p^A_M:T^AM\to M$ for the bundle projection. Since
$S M\subset TTM\car_M TTM$ is defined by (\ref{e11}),
$T^AS M\subset T^ATTM\car_{T^AM}T^ATM$
is the set of all pairs $(X,Y)$ satisfying
\bEq\label{e15}
T^Ap^T_{TM} X= T^ATp^T_M Y\,,\quad T^Ap^T_{TM} Y= T^ATp^T_M X\,.
\eEq
On the other hand, $S T^AM\subset TTT^AM\car_{T^A M} TTT^AM$ is
characterized by
\bEq\label{e16}
    p^T_{TT^AM} X= Tp^T_{T^AM} Y\,,\quad p^T_{TT^AM} Y= Tp^T_{T^AM} X\,.
\eEq

We have $T^A\sig_M:T^AS M \to T^A TM$,
$\kappa^A_{TM}:T^ATTM\to TT^ATM$
 and
$T\kappa^A_M:TT^A TM\to TTT^AM$. For technical reasons, we postpone
the proof of the following assertion to Section 5.

\bPr\label{Pr2}
The map $T\kappa^A_M\com \kappa^A_{TM}:T^A TTM\to TTT^AM$
induces a diffeomorphism $K^A_M:T^AS M\to S T^AM$ and the following
diagram commutes
\bEq\label{e17}
\begin{CD}
T^AS M @>K^A_M>> S T^AM
\\
@V{T^A\sig_{M}}VV @VV{\sig_{T^AM}}V
\\
T^ATM @>\kappa^A_M>> TT^AM
\end{CD}
\eEq
\ePr

Now we first show how (\ref{e17}) implies
that the flow prolongation $\mc{T}^A$ of classical vector
fields $X,Y:M\to TM$ preserves the bracket. We have
$(TY\com X, TX\com Y):M\to S M$ and
\bEq\label{e18}
[X,Y]=\sig_M\com(TY\com X, TX\com Y)\,.
\eEq
Then $T^A(TY\com X, TX\com Y):T^AM\to T^AS M$. Adding
$K^A_M$ we obtain
$T\kappa^A_M\com \kappa^A_{TM}\com T^ATY\com T^AX
        = T\kappa^A_M\com TT^AY\com \kappa^A_M\com T^A X
        = T\mc{T}^A Y\com \mc{T}^A X$
and the same for $TX\com Y$.
So in (\ref{e17}) we clockwise obtain
$[\mc{T}^AX,\mc{T}^A Y]$.
Counterclockwise, we first get $T^A[X,Y]$ and then
$\mc{T}^A[X,Y]$.

Consider now the case of $\F$. According to the general fact that the
homomorphisms of Weil algebras extend to the functional case,
(\ref{e17}) yields a commutative diagram
\bEq\label{e19}
\begin{CD}
T^AS \F @>K^A_{\F}>> S T^A\F
\\
@V{T^A\sig_{\F}}VV @VV{\sig_{T^A\F}}V
\\
T^AT\F @>\kappa^A_{\F}>> TT^A\F
\end{CD}
\eEq

\noindent
For two vector fields $X,Y$ on $\F$, we first construct
\beq
    (TY\com X, TX\com Y):\F\to S \F\,.
\eeq
Then we deduce (\ref{e10}) in the same way as in the manifold
case. This proves Proposition \ref{Pr1}.

\section{Some Weilian lemmas}
\setcounter{equation}{0}
\setcounter{theorem}{0}

The elements of $A=T^A\Rn$ are of the form $j^Ag$, $g:\Rn^k\to \Rn$.
For a vector space $V$, the map $V\car A\to T^AV$,
$(v,j^Ag)\mto j^A(gv)$ is bilinear and defines an identification
$T^AV= V\ten A$.
If $W$ is another vector space and $f:V\to W$ is a linear map, then
$T^Af : T^AV \to T^AW$ is of the form
\bEq\label{e20}
    T^Af=f\ten \id_A:V\ten A\to W\ten A\,,
\eEq
\cite{KolMicSlo93}. Further, let $\mu:A\to B$ be an algebra
homomorphism. Then the induced natural transformation
$\mu_V:T^AV\to T^BV$
is of the form
\bEq\label{e21}
\mu_V=\id_V\ten \mu:  V\ten A\to V\ten B\,.
\eEq
This follows from the fact that $V$ is isomorphic to $\Rn^n$
and we have a product preserving functor.

In particular, if $C$ is another Weil algebra, then (\ref{e20})
implies that the natural transformation $T^C\mu_M:T^CT^AM\to T^CT^BM$
corresponds to the algebra homomorphism
\bEq\label{e22}
\id_C\ten \mu: C\ten A\to C\ten B\,.
\eEq
Further, the maps $\mu_{T^CM}:T^AT^CM\to T^BT^CM$ form a natural
transformation $T^AT^C\to T^BT^C$ that corresponds
to the algebra homomorphism
\bEq\label{e23}
\mu\ten \id_C: A\ten C\to B\ten C\,.
\eEq

The trivial bundle functor on $\mc{M}f$ transforming every
manifold $M$ into $\id_M:M\to M$ and every smooth map $f$
into $(f,f)$ corresponds to the trivial Weil algebra $\Rn$.
The natural transformation $p^A_M:T^AM\to M$ is determined
by the canonical "real part projection"
$\rho_A:A=\Rn\car N_A\to \Rn$.
So $T^Bp^A_M:T^BT^AM\to T^BM$ corresponds to the canonical map
\bEq\label{e24}
\id_B\ten\rho_A: B\ten A\to B\ten \Rn=B\,.
\eEq

Write $\kappa^{A,B}:A\ten B\to B\ten A$ for the exchange map.
This defines the exchange natural transformation
$\kappa^{A,B}_M:T^AT^B M\to T^BT^A M$. By (\ref{e23}),
$\kappa^{A,B}_{T^C M}:T^AT^BT^C M\to T^BT^AT^C M$
corresponds to the exchange $A\ten B\ten C\to B\ten A\ten C$.
By (\ref{e22}), $T^B\kappa^{A,C}_M:T^BT^AT^CM\to T^BT^CT^A M$
corresponds to the exchange $B\ten A\ten C\to B\ten C\ten A$.

\bLm
The following diagram commutes
\bEq\label{e25}
\begin{CD}
T^AT^BT^CM @>{T^B\kappa^{A,C}_M\com\kappa^{A,B}_{T^CM}}>> T^BT^CT^AM
\\
@V{T^A p^{B}_{T^CM}}VV @VV{p^B_{T^CT^AM}}V
\\
T^AT^CM @>\kappa^{A,C}_M>> T^CT^AM
\end{CD}
\eEq
\eLm

\bPf
At the algebra level, we have a commutative diagram

\medskip
\begin{center}
\setlength{\unitlength}{1mm}
\begin{picture}(70,25)
\put(0,23){\makebox(0,0)[cc]{$A\ten B\ten C$}}
\put(35,23){\makebox(0,0)[cc]{$B\ten A\ten C$}}
\put(70,23){\makebox(0,0)[cc]{$B\ten C\ten A$}}
\put(0,3){\makebox(0,0)[cc]{$A\ten C$}}
\put(70,3){\makebox(0,0)[cc]{$C\ten A$}}
\put(11,23){\vector(1,0){13}}
\put(46,23){\vector(1,0){13}}
\put(0,19){\vector(0,-1){12}}
\put(70,19){\vector(0,-1){12}}
\put(15,3){\vector(1,0){40}}
\end{picture}
\end{center}
\ePf

\vglue-5mm

Now we are in position to prove Proposition \ref{Pr2}.
Comparing our general case with the situation in Section 4,
we see $\kappa^{A,\mathbb{D}}=\kappa^A$ and
$p^{\mathbb{D}}_M=p^T_M$. So if we put $B=\mathbb{D}=C$
into (\ref{e25}), we obtain
\bEq\label{e26}
    p^T_{TT^A M}\com T\kappa^A_M\com \kappa^A_{TM}
        = \kappa^A_M\com T^Ap^T_{TM}\,.
\eEq
Every $X,Y\in T^AS M$ satisfy (\ref{e15}).
The naturality of $\kappa^A$ on $p^T_M:TM\to M$ yields
\bEq\label{e27}
    \kappa^A_M\com T^ATp^T_M
        = TT^Ap^T_M\com \kappa^A_{TM}
\eEq
and the standard relation $p^T_{T^AM}\com \kappa^A_M= T^Ap^T_M$ implies
\bEq\label{e28}
    Tp^T_{T^A M}\com T\kappa^A_M
        = TT^Ap^T_M\,.
\eEq
Hence we have $(p^T_{TT^AM}\com T\kappa^A_M\com \kappa^A_{TM})(X)
=\kappa^A_M(T^Ap^T_{TM}(X))=\kappa^A_M(T^ATp^T_M(Y))=
(TT^Ap^T_M\com\kappa^A_{TM})(Y)
=(Tp^T_{T^AM}\com T\kappa^A_M\com \kappa^A_{TM})(Y)$.
Thus, $(T\kappa^A_M\com \kappa^A_{TM})(X)$ and
$(T\kappa^A_M\com \kappa^A_{TM})(Y)$ satisfy (\ref{e16}),
so that $K^A_M$ maps $T^AS M$ into $S T^A M$.
In the case $M= \Rn^m$, we have $S \Rn^m=\overset{5}{\car}\Rn^m$
and $T^A\Rn^m=A^m$, so that $T^AS \Rn^m=\overset{5}{\car}A^m$
and $S T^A\Rn^m=\overset{5}{\car}A^m$.
In this situation, $K^A_{\Rn^m}$ is the identity of
$\overset{5}{\car}A^m$.
Moreover,
by (\ref{e14}) $\sig_{\Rn^m}$ is determined by the difference
of the fourth and fifth components. Taking into account that
the vector addition in $A$ is the $T^A$-prolongation of the
addition of reals, we deduce that the diagram (\ref{e17}) commutes.

\section{The jet-like case}
\setcounter{equation}{0}
\setcounter{theorem}{0}

Every fiber product preserving bundle functor $G$ on $\mc{F}\mc{M}_m$
is of the form $G=(A,H,t)$ where $A$ is a Weil algebra, $H:G^r_m\to
\Aut A$ is a group homomorphism and $t:\mathbb D^r_m\to A$ is
an equivariant algebra homomorphism, \cite{KolMik99}.
For every manifold $N$, the natural transformations corresponding
to $\Aut A$ determine an action
$H_N$ of $G^r_m$ on $T^AN$. So we can construct
the associated bundle $P^rM[T^AN,H_N]$, where
$P^rM\subset T^r_m M$ is the $r$-th order frame bundle of $M$.
For a fibered manifold $\pi:E\to M$, we define $GE$ as a subset of
$P^rM[T^AE,H_E]$ characterized by
\bEq\label{e29}
    GE=\{\{u,Z\}, t_Mu=T^A\pi(Z)\}\,,\quad
   u\in P^rM, Z\in T^AE\,.
\eEq
For an $\mc{F}\mc{M}_m$-morphism $f:E\to\bar E$ over a local
diffeomorphism $\underline{f}:M\to\bar M$, we have the induced
principal bundle morphism
$P^r\underline{f}:P^rM\to P^r\bar M$ and an $G^r_m$-equivariant
map $T^Af:T^AE\to T^A\bar E$. So we can construct
$P^r\underline{f}[T^Af]:P^rM[T^AE]\to P^r\bar M[T^A\bar E]$
and we define
\bEq\label{e30}
Gf= P^r\underline{f}[T^Af]|GE\,.
\eEq
In the product case $E=\Rn^m\car Q$, we have
$GE= \Rn^m\car T^AQ$, \cite{KolMik99}.

This construction extends directly to $\F$. By (\ref{e8}),
each element of $\Aut A$ determines an $F$-smooth
isomorphism $T^A\F\to T^A\F$. So we have an action
$H_{\F}$ of $G^r_m$ on $T^A\F$ and we can construct
the $F$-smooth associated bundle
\bEq\label{e31}
    P^rM[T^A\F,H_{\F}]\,.
\eEq
Then we define $G\F$ as the subset of (\ref{e31}) characterized by
\bEq\label{e32}
    G\F=\{\{u,Z\}, t_Mu= T^Ap(Z)\}\,,\quad u\in P^rM, Z\in T^A\F\,.
\eEq
Write $\mc{F}\mc{M}^I_m=\mc{F}\mc{M}^I\cap \mc{F}\mc{M}_m$.
For $(f_1,f_2)\in\mc{F}\mc{M}^I_m\car_{\mc{B}} \mc{F}\mc{M}_m$
with the common base map $\underline{f}$, we define
\bEq\label{e33}
    G\mc{F}(f_1,f_2) = P^r\underline{f}[T^A\mc{F}(f_1,f_2)]|G\F\,.
\eEq
Hence $G\mc{F}$ is a functor on
$\mc{F}\mc{M}^I_m\car_{\mc{B}} \mc{F}\mc{M}_m$
with values in $\mc{S}\mc{B}$.

In the product case $E_1= \Rn^m\car Q_1$, $E_2= \Rn^m\car Q_2$, we
have
\bEq\label{e34}
    G\F= \Rn^m\car C^\infty(Q_1,T^AQ_2)\,.
\eEq
This shows that for
$J^r= (\mathbb{D\,}{}^r_m, \id_{G^r_m}, \id_{\mathbb{D\,}{}^r_m})$
we obtain $J^r\F$ constructed by means of the fiber $r$-jets
in \cite{CabKol95}.

\section{Vector fields in the jet-like case}
\setcounter{equation}{0}
\setcounter{theorem}{0}

In the manifold case, \cite{Kol03}, if we have a principal bundle
$P(M,C)$ with structure group $C$ and a left $C$-space $S$,
a right-invariant vector field $\varphi$ on $P$ and a left-invariant
vector field $\psi$ on $S$, the product vector field
$(\varphi,\psi)$ on $P\car S$ is projectable to a vector field
$\{\varphi,\psi\}$ on the associated bundle $P[S]$.
In particular, if $\eta$ is a projectable vector field
on $E\to M$ over a vector field $\xi$ on $M$,
then the flow prolongation $\mc{P}^r\xi$ is right-invariant
on $P^rM$ and $\mc{T}^A\eta$ is left-invariant on $T^A E$.
In \cite{Kol03} we deduced that the flow prolongation $\mc{G}\eta$
of $\eta$
coincides with the restriction of $\{\mc{P}^r\xi,\mc{T}^A\eta\}$
to $GE\subset P^rM[T^AE]$.

In the functional case, consider a vector field
$X:\F\to T\F$ over $\xi:M\to TM$. Then (\ref{e8})
implies that the field prolongation $\mc{T}^AX$ is $H_{\F}$-invariant.
Hence we have the vector field $\{\mc{P}^r\xi,\mc{T}^AX\}$
on $P^rM[T^A\F]$ and we define the field prolongation $\mc{G}X$
of $X$ by
\bEq\label{e35}
\mc{G}X=\{\mc{P}^r\xi,\mc{T}^AX\}|G\F\,.
\eEq
This is a vector field $G\F\to TG\F$ over $\xi$. For
two vector fields $X_i$ on $\F$ over $\xi_i$, $i=1,2$, we have
by the basic properties of the strong difference
\beq
    [\mc{G}X_1,\mc{G}X_2]
        = \{[\mc{P}^r\xi_1,\mc{P}^r\xi_2],[\mc{T}^AX_1,\mc{T}^AX_2]\}\,.
\eeq
Hence Proposition \ref{Pr1} yields

\bPr\label{Pr3}
We have
\beq
    [\mc{G}X_1,\mc{G} X_2] = \mc{G}[X_1,X_2]\,.\quad \END
\eeq
\ePr

At the end we remark that the third author, \cite{Kol03},
constructed a map
\beq
    \mu_E^G:J^rTM\car_{GTM}GTE\to TGE
\eeq
with the property that for every projectable vector field
$\eta$ on $E$ over $\xi$ on $M$
\beq
    \mc{G}\eta=\mu^G_E\com(j^r\xi\car_M G\eta)\,,
\eeq
where $j^r\xi:M\to J^rTM$ is the $r$-th jet prolongation
of the section $\xi:M\to TM$ and $G\eta:GE\to GTE$ is the induced
morphism. Analyzing this construction, one realizes that each step
can be extended to our functional case.
In other words, one can introduce in the same way
an $F$-smooth morphism
\beq
    \mu^G_{\F}:J^rTM\car_{GTM}G\F\to TG\F
\eeq
with the property
\beq
    \mc{G}X=  \mu^G_{\F}\com (j^r\xi\car_M GX)
\eeq
for every vector field $X$ on $\F$ with underlying vector field
$\xi$ on $M$.



\end{document}